\input amstex
\documentstyle{amsppt}

\def\R{\Bbb R}
\def\Z{\Bbb Z}

\def\Cplx{\Bbb C}

\def\hy{\hbox{\it -}}
\def\spec{\text{spec}}

\def\Ad{\text{Ad}}
\def\ad{\text{ad}}
\def\and{\text{and}}
\def\span{\text{span}}
\def\Aut{\text{Aut}}
\def\Der{\text{Der}}
\def\OA{\text{OA}}
\def\IA{\text{IA}}
\def\ID{\text{ID}}
\def\AIA{\text{AIA}}
\def\AID{\text{AID}}
\def\GW{{\text{GW}}}
\def\F{{F_\tau}}
\def\Tau{\Cal T}
\def\grad{\text{grad}}
\def\H{\Cal H}
\def\CZ{\Cal Z}

\def\SL{\text{SL}}

\def\OD{\text{OD}}
\def\isom{=}

\def\gp{G}
\def\Center{Z(\gp)}
\def\dgp#1{\gp^{(#1)}}
\def\GG{\Gamma}
\def\G#1{\GG_{#1}}
\def\Nilmfld{\GG\backslash\gp}
\def\nilmfld#1{\G{#1}\backslash\gp}

\def\alg{\frak g}
\def\dalg#1{\alg^{(#1)}}
\def\center{\frak z}

\def\nn{\frak N}
\def\Char{\frak C}

\def\qg{\bar{g}}
\def\qgp{\bar\gp}
\def\qGG{\bar\GG}
\def\qG#1{\qGG_{#1}}
\def\qNilmfld{\qGG\backslash\qgp}
\def\qnilmfld#1{\qG{#1}\backslash\qgp}

\def\qalg{\bar\alg}

\def\a{\alpha}
\def\z{\zeta}

\def\hg{\hat{g}}
\def\hgp{\hat{\gp}}
\def\hGG{\hat{\GG}}
\def\hG#1{\hGG_{#1}}
\def\hNilmfld{\hGG\backslash\hgp}
\def\hnilmfld#1{\hG{#1}\backslash\hgp}
\def\hfnsp#1{L^2(\hnilmfld{#1})}
\def\halg{\hat\alg}
\def\hX{\hat{X}}
\def\hZ{\hat{Z}}

\define\Rep{{\raise 1.5pt\hbox{$\rho$}_\GG}}
\define\rep#1{{\raise 1.5pt\hbox{$\rho$}_{\G{#1}} }}
\define\hRep{{ \raise 1.5pt\hbox{$\rho$}_\hGG }}
\define\hrep#1{{ \raise 1.5pt\hbox{$\rho$}_{\hG{#1}} }}

\topmatter

\NoBlackBoxes
\TagsAsMath

\pagewidth{4.7 in}
\pageheight{7.5 in}
\magnification=1100

\title 
Continuous Families of Riemannian manifolds, 
isospectral on functions but not on 1-forms.
\endtitle

\rightheadtext
{Continuous Families of Isospectral Riemannian Manifolds.}

\author
Ruth Gornet
\endauthor

\address
Ruth Gornet; Texas Tech University; Department of Mathematics; 
Lubbock, Texas \ 79409-1042; e-mail: gornet\@math.ttu.edu
\indent February 1997
\endaddress

\keywords
Isospectral manifolds; Laplacian on forms; 
Nilmanifolds; Isospectral deformations; Marked length spectrum.
\newline
Research at Texas Tech University supported in part by a 
State of Texas, Advanced Research Program Grant
\endkeywords


\abstract
The purpose of this paper is to present the
first continuous families of Riemannian manifolds
isospectral on functions but not on 1-forms,
and simultaneously, the first continuous
families of Riemannian manifolds with the same marked length
spectrum but not the same 1-form spectrum. 
The examples presented here are Riemannian three-step nilmanifolds
and thus provide a 
counterexample to the Ouyang-Pesce Conjecture for higher-step 
nilmanifolds.
Ouyang and Pesce independently showed that all isospectral
deformations of two-step nilmanifolds must 
arise from the Gordon-Wilson method for constructing isospectral nilmanifolds.
In particular, all continuous families of Riemannian two-step
nilmanifolds that are isospectral on functions must
also be isospectral on $p$-forms for all $p.$
They conjectured that all isospectral deformations of
nilmanifolds must arise in this manner.  
These examples arise from a general method for constructing
isospectral Riemannian nilmanifolds
previously introduced by the author.
\endabstract

\endtopmatter

\document

\subheading{Section 1. Introduction}
\medskip

The {\bf spectrum} of a closed Riemannian manifold $(M,g)$, 
denoted $\spec(M,g)$, is the collection of eigenvalues with 
multiplicities of the associated Laplace--Beltrami operator 
acting on smooth functions.
Two Riemannian manifolds $(M,g)$ and $(M',g')$ are said to be
{\bf isospectral} if $\spec(M,g)=\spec(M',g').$
The Laplace--Beltrami operator may be extended to act on smooth 
$p$-forms by $\Delta = d \delta + \delta d,$  where $\delta$ 
is the metric adjoint of $d.$  We denote its eigenvalue spectrum
by $p\hy\spec(M,g).$  

This paper addresses the question of spectral rigidity:
``If we continuously deform a manifold while
keeping the Laplace spectrum constant, must the resulting manifolds
be isometric?'' That is, can we determine the geometry
of a manifold based on nearby behavior of the spectrum?

Sometimes the answer is ``Yes.'' Croke and Sharafutdinov \cite{CS},
generalizing work of Guillemin and Khazdan \cite{GuK} for surfaces,
have recently showed spectral rigidity of negatively curved manifolds.
Flat tori are known to be spectrally rigid \cite{Wp}, 
while Heisenberg manifolds \cite{P1} are spectrally rigid 
within the family of left invariant metrics.  

Nilmanifolds and nilpotent Lie groups have played a vital
role in showing that the general answer is ``No.''
In 1984, Gordon and Wilson \cite{GW1} constructed 
the first nontrivial isospectral deformations,
continuous families of isospectral, nonisometric Riemannian nilmanifolds.
The Gordon--Wilson method has been modified to produce conformally 
equivalent isospectral deformations of Riemannian 
two-step nilspaces \cite{DG2},
\cite{BG}, and nontrivial isospectral deformations 
arbitrarily close to Heisenberg manifolds \cite{Sch2}.
Using new techniques, Gordon and Wilson \cite{GW3}
have recently used nilmanifolds to construct continuous families
of isospectral (Dirichlet and Neumann)
manifolds with boundary that are nonlocally isometric.
These are the only known examples of nontrivial isospectral
deformations.

While two-step nilmanifolds are not spectrally rigid,
they do satisfy a ``modified'' rigidity.
In 1992, Ouyang \cite{O} and Pesce \cite{P3},
independently showed that isospectral deformations of two-step nilmanifolds
must arise from the Gordon--Wilson method. 
Ouyang and Pesce conjectured, 
``{\it isospectral deformations of
nilmanifolds must arise from the Gordon--Wilson method.}''
It was expected that a combination of their different techniques
would prove the result on higher-step nilmanifolds.

In 1985, Sunada \cite{Su} gave an elegant method for constructing pairs of
isospectral manifolds that can be explained solely in representation
theoretic terms. (See \cite{B2} for a discussion.)
His method was later generalized \cite{DG2} to include
the Gordon--Wilson deformations. Sunada's method and
its generalizations produce manifolds that are {\it strongly} isospectral,
that is, they must share the same spectrum for all natural, strongly
elliptic operators. In particular, isospectral deformations
arising from the Gordon--Wilson method
must share the same $p$-form spectrum for all $p.$
Conversely,  Schueth \cite{Sch1} showed that if
one uses representation theoretic techniques alone
to obtain continuous families of isospectral nilmanifolds,
then the deformation must arise from the Gordon--Wilson method.

As most constructions for producing families of 
isospectral manifolds can by explained by Sunada's method or
its generalizations,
examples of Riemannian manifolds that are isospectral on functions
but not on forms are sparse. 
Gordon \cite{G2} has constructed
pairs of Heisenberg manifolds that are isospectral on functions but 
not on 1-forms, and the author \cite{Gt3,4} has constructed pairs of three-step
Riemannian nilmanifolds with this property.  
The only other known examples of manifolds 
that are isospectral on functions but not on $p$-forms for all $p$
are lens spaces. 
For any choice of $P\in\Z^+,$
Ikeda \cite{I2} has constructed examples of pairs isospectral lens
spaces that are isospectral on $p$-forms for $p=0,1,\dots, P$
but not isospectral on $(P+1)$-forms. 
Recently, Pesce \cite{P4} has been able to explain the Ikeda examples
in a Sunada-like setting, which requires a genericity assumption
excluding nilmanifolds. 
Nilmanifolds, therefore, now account for all known isospectral
manifolds that do not fall under a Sunada-like setting.

\smallskip
The main result of this paper is the following.
\smallskip

\flushpar{\bf Theorem: 
There exist continuous families of Riemannian manifolds
that are isospectral on functions but not on 1-forms.}
\smallskip

The examples presented here are three-step nilmanifolds,
and thus provide a counterexample to the Ouyang--Pesce conjecture.
The new deformations arise from a general method 
for constructing isospectral Riemannian nilmanifolds introduced in \cite{Gt4},
which combines techniques from
representation theory and Riemannian geometry.

All known examples of strongly isospectral nilmanifolds
arise from the Gordon--Wilson method, leading to the following.
\smallskip

\flushpar {\bf Conjecture: Continuous families of
nilmanifolds that are isospectral on functions and 
on 1-forms must arise from the Gordon--Wilson method.} 
\smallskip

Two Riemannian manifolds $(M_1,g_1)$ and $(M_2,g_2)$ have the 
same {\bf marked length spectrum} if there exists an isomorphism
between the fundamental groups of $M_1$ and $M_2$ such that
corresponding free homotopy classes of loops contain smoothly
closed geodesics of the same length. (See Section 5 for more details.)

The relationship between the Laplace spectrum and lengths of closed
geodesics arises from the study of the wave equation (see \cite{DGu}),
and in the case of compact, hyperbolic manifolds, from the Selberg
Trace Formula (see \cite{C, Chapter XI}). Colin de Verdi\`ere \cite{CdV}
has shown that generically, the Laplace spectrum determines the lengths
of closed geodesics. All known examples of manifolds that are isospectral
have the same lengths of closed geodesics.

Croke \cite{Cr} and Otal \cite{Ot12} independently showed that if a 
pair of compact surfaces with negative curvature have the same marked
length spectrum, they are necessarily isometric. This is due to
the fact that in certain cases, the marked length spectrum
and the geodesic flow are (roughly) equivalent notions.
Gordon, Mao, and Schueth \cite{GM1,2}, \cite{GMS} have shown that
generically, two-step nilmanifolds with conjugate geodesic flows
must be isometric. 

Isospectral Riemannian manifolds constructed using the Gordon--Wilson
method necessarily have the same marked length spectrum, while
Eberlein \cite{E} proved the converse in the case of two-step nilmanifolds. 
Thus the Gordon--Wilson method produces continuous families of 
strongly isospectral Riemannian manifolds with the same marked length spectrum.
In contrast, we prove the following.

\flushpar{\bf Theorem: 
There exist continuous families of Riemannian manifolds
with the same marked length spectrum but not the same spectrum on 1-forms.}
\smallskip

Discrete families of manifolds with this property were presented in \cite{Gt4},
to which the reader is referred for an introduction to the  
relationships among the marked length spectrum, the length spectrum,
the Laplace spectrum on functions, and the Laplace spectrum on forms on
discrete families of Riemannian manifolds.
The only other known examples of manifolds that have the same marked length spectrum
but not the same spectrum on 1-forms are the Zoll spheres and standard spheres,
which are not even isospectral on functions.

This paper is organized as follows.
We briefly review background material and the Gordon--Wilson method
in Section 2. The techniques used to construct the new examples
are discussed in Section 3. One way to show that isospectral nilmanifolds
cannot arise from the Gordon--Wilson method is by showing 
that no diffeomorphism between the manifolds of the correct 
type can exist. This approach is discussed in Section 3,
as it lends insight into the types of Riemannian nilmanifolds
that can support isospectral deformations that are not strongly 
isospectral.
The new examples are presented in Section 4,
while the marked length spectrum and the 1-form spectrum are the
focus of Sections 5 and 6, respectively.

\subheading{Acknowledgments}
The author wishes to thank the joint NSF-CNRS program in Spectral Geometry
for their support of the workshop at Dartmouth in February 1996, where
much of this work was carried out. The author also wishes to thank
Ed Wilson and David L. Webb for insightful conversations.

\bigskip

\subheading{Section 2. Background and Notation}
\medskip

\subheading{2.1 Riemannian Nilmanifolds}
Let $\gp$ be a simply connected Lie group, and let 
$\GG$ be a cocompact, discrete subgroup of $\gp.$ A Riemannian metric 
$g$ is {\bf left invariant} if the left translations of $\gp$ are isometries.   
The left invariant metric $g$ descends to a Riemannian metric on $\Nilmfld,$ 
also denoted by $g.$  Note that a left invariant metric is 
determined by specifying an orthonormal basis
$\{E_1, \cdots, E_n\}$  of the Lie algebra $\alg$ of $\gp.$
As $\gp$ is unimodular, the Laplace--Beltrami
operator $\Delta$ of $(\Nilmfld, g)$ is just 
$$\Delta = -\sum_{i=1}^n {E_i}^2.$$ 

For a Lie algebra $\alg$, denote by $\dalg1$  
the derived algebra $[\alg,\alg]$ of $\alg.$ 
Inductively, define $\dalg{k+1}=[\alg,\dalg{k}].$
A Lie algebra
$\alg$ is {\bf k-step nilpotent} if $\dalg{k-1}$ is nontrivial and central.
A Lie group $\gp$ is {$k$-step nilpotent} if its Lie algebra is.
Let $\GG$ be a cocompact,
discrete subgroup of a simply connected nilpotent Lie group 
$\gp$ with left invariant metric $g.$   
The locally homogeneous space 
$(\Nilmfld,g)$ is called a {\bf Riemannian nilmanifold.}

For more details about nilmanifolds and nilpotent Lie groups see 
\cite{CG}.

\subheading{2.2 Isometries of Nilmanifolds} 
We write $(M,g)\isom (M',g')$ if the Riemannian manifolds
$(M,g)$ and $(M',g'),$ are isometric.
Let $(\Nilmfld, g)$ be a Riemannian nilmanifold.
Any automorphism $\Phi\in\Aut(\gp)$ 
induces an isometry from 
$(\Nilmfld,\Phi^*g)$ to $(\Phi(\GG)\backslash\gp,g).$
Any inner automorphism $I_x\in\IA(\gp)$
induces an isometry from $(\Nilmfld,g)$ to $(\Nilmfld, I_x^*g)$
via the mapping $R_x:\Nilmfld\to\Nilmfld$ defined by $R_x(\GG y)=\GG yx.$
Let $\OA(\gp,g)=\{\Upsilon\in\Aut(\gp): \Upsilon^*g=g\},$
the orthogonal automorphisms of $\gp$ with respect to $g.$
If $\Upsilon$ is in $\OA(\gp,g),$
then $(\Nilmfld,g) \isom (\Nilmfld,\Upsilon^*g).$

Recall that for any Lie group $\gp$ with
Lie algebra $\alg$, $\Aut(\alg)$ is a Lie group with Lie
algebra $\Der(\alg).$ 
We denote by $\exp:\Der(\alg)\to\Aut(\alg)$ 
the Lie algebra exponential,
which reduces to the matrix exponential by looking at matrix
representations.

\definition{2.3 Notation}
For every $D\in\Der(\alg),$
we associate an operator $D^*$ that takes the two-form $\omega$
to the two-form defined by $D^*\omega(X,Y)=\omega(DX,Y)+\omega(X,DY).$
\enddefinition

One generally
denotes by $\Phi^*\omega$ the two-form defined by 
$\Phi^*\omega(X,Y)=\omega(\Phi_*(X),\Phi_*(Y)),$
where $\Phi_*\in\Aut(\alg).$
As $D\in\Der(\alg)$ is the derivative of a curve $\Phi_{s*}$
in $\Aut(\alg),$ $D^*\omega$ is the derivative of 
$\Phi_s^*\omega$ in the obvious sense.

With this notation, $\OA(\gp,g)$ is a Lie group with Lie algebra
$\OD(\gp,g)=\{D\in\Der(\alg): D^*g=0\},$
the skew-symmetric derivations of $\alg.$

\definition{2.4 Almost Inner Automorphisms \cite{GW1}, \cite{G1}} 
Let $\gp$ be a Lie group with Lie algebra $\alg$
and cocompact, discrete subgroup $\GG.$
A Lie group automorphism $\Psi$ of $\gp$ is {\bf almost inner}
if for all $x$ in $\gp,$ the elements
$x$ and $\Psi(x)$ are conjugate in $\gp.$
The automorphism $\Psi$ of $\gp$ is ${\bold\GG}$-{\bf almost inner} if 
for all $\gamma$ in $\GG,$ the elements
$\gamma$ and $\Psi(\gamma)$ are conjugate in $\gp.$ 
A Lie algebra isomorphism $\Psi_*$
is {\bf almost inner} (resp., ${\bold\GG}$-{\bf almost inner})
if for every $X$ in $\alg,$ (resp.,  $X$ in $\log\GG$),
$\Psi_*(X)=Ad(x)(X)$ for some $x$ in $\alg.$ 
Here $\Ad$ denotes the adjoint representation of $\gp$ on $\alg.$ 
A derivation $D$ is inner if there exists $Y$ in $\alg$ such
that $D=\ad(Y).$  A derivation $D$ is {\bf almost inner}, (resp.,
${\bold\GG}$-{\bf almost inner}) if for every $X$ in $\alg$ (resp.,
$X$ in $\log\GG$) there exists $Y$ in $\alg$ such that $D(X)=\ad(Y)(X).$
\enddefinition

Denote by $\IA(\gp)$, (resp., $\AIA(\gp), \AIA(\gp;\GG)$ )
the set of inner 
(resp., almost inner, $\GG$-almost inner) automorphisms of $\gp.$
Note that 
$\IA(\gp)\subset\AIA(\gp)\subset\AIA(\gp;\GG),$
and there exist examples \cite{Gt1}, \cite{GW1} 
demonstrating that each of these containments can be strict.
Define $\IA(\alg),$ $\AIA(\alg),$ $\AIA(\alg;\GG),$
$\ID(\alg)$, $\AID(\alg),$ and $\AID(\alg; \GG)$ analogously.

If $\gp$ is nilpotent, then $\Psi\in\Aut(\gp)$ is 
inner (resp., almost inner, $\GG$-almost inner) 
if and only if $\Psi_*\in\Aut(\alg)$ is
inner (resp., almost inner, $\GG$-almost inner). 
Also, $\IA(\alg),$ $\AIA(\alg)$, and $\AIA(\alg;\GG)$
are connected Lie subgroups of $\Aut(\alg),$
with Lie algebras $\ID(\alg),$ $\AID(\alg),$ and $\AID(\alg; \GG),$ 
respectively. 

\proclaim{2.5 Theorem \cite{GW1}, \cite{G2}}
Let $\GG$ be a cocompact, discrete subgroup of a simply connected, nilpotent 
Lie group $\gp$ with left invariant metric $g.$ 
\roster
\item If $\Psi$ is a $\GG$-almost inner automorphism,
then $p\hy\spec(\Nilmfld,g) = p\hy\spec(\Psi(\GG) \backslash\gp,g)$
for $p=0,1,\cdots,$ $\dim(\gp).$
\item If $\Psi_s$ is a continuous family of $\GG$-almost
inner automorphisms that are not inner automorphisms,
then $(\Psi_s(\GG)\backslash\gp,g)=(\Nilmfld,\Psi_s^*g)$ 
is a nontrivial, isospectral deformation.
\item There exist two-step nilpotent Lie groups supporting
continuous families of $\GG$-almost inner automorphisms
that are not inner automorphisms.
\endroster
\endproclaim

Conversely, Ouyang and Pesce independently proved
the following.

\proclaim{2.6 Theorem \cite{O},\cite{OP},\cite{P2}}
Let $(\Nilmfld, g_s)$ be a continuous family of Riemannian
two-step nilmanifolds that are isospectral on functions.
Then there exists a continuous family $\Psi_s$ of 
$\GG$-almost inner automorphisms
with $\Psi_0=I$ such that $g_s=\Psi_s^*g_0.$
\endproclaim

\subheading{2.7 Notation and Remarks}
Define 
$$\GW(\gp,g; \GG)=\OA(\gp,g)\circ\AIA(\gp; \GG)=\{\Upsilon\circ\Psi:
\Upsilon\in\OA(\gp,g), \Psi\in\AIA(\gp;\GG)\},$$
 
$$\GW(\gp,g)=\OA(\gp,g)\circ\AIA(\gp)=\{\Upsilon\circ\Psi:
\Upsilon\in\OA(\gp,g), \Psi\in\AIA(\gp)\}.$$

The group $\AIA(\gp)$ is a normal subgroup 
of $\Aut(\gp),$ so $\GW(\gp,g)$ is a Lie group with Lie algebra
$\OD(\alg,g)\oplus \AID(\alg).$
In particular, if $D\in\OD(\alg,g)\oplus \AID(\alg),$
then the Lie group automorphism $\Phi$ defined at the
Lie algebra level by $\Phi_*=\exp(D)$ is in $\GW(\gp,g).$
A priori, the set $\GW(\gp,g;\GG)$ is not a subgroup of $\Aut(\gp)$,
as there is no reason to assume that $\AIA(\gp; \GG)$ is a normal
subgroup of $\Aut(\gp),$ although this proves to be the case in
all known examples.

For the nilpotent Lie groups studied in Section 4, 
$\AIA(\gp)=\AIA(\gp; \GG)$ for any choice of cocompact, discrete
subgroup $\GG$ of $\gp.$

\proclaim{2.8 Proposition \cite{GW1}}
Let $\gp$ be a simply connected, nilpotent Lie group with
cocompact discrete subgroup $\GG$ and left invariant metric $g.$
Let $\Phi_s$ be a continuous family of automorphisms of $\gp$
such that $\Phi_0$ is the identity.  
\roster
\item The continuous family $(\Nilmfld, \Phi_s^*g)$ is trivial
if and  only if $\Phi_s\in\OA(\gp,g)\circ\IA(\gp)$
for all $s.$
\item The continuous family $(\Nilmfld, \Phi_s^*g)$ 
is Gordon--Wilson if and only if there exists $\Upsilon$ in $\OA(G,g)$ 
such that 
for all $s,$ $\Phi_s\circ\Upsilon\in\GW(\gp,g;\GG).$
\endroster
\endproclaim

Note that by (2.7), if $\GW(\gp; \GG)$ is a group,
then we may assume $\Upsilon=I$ in
the statement of Proposition 2.8.

\demo{Comments on Proof}
This follows directly from \cite{GW1, Corollary 5.3},
using  $\Cal D=\{\delta\in\Aut(\gp): \delta(\GG)=\GG\}$ 
discrete and $\Phi_0$ the identity.
Note that the group defined by Gordon--Wilson as
$\text{K}_g$ is denoted here by $\OA(\gp,g).$ 
\qed \enddemo

\subheading{2.9 Quotient Nilmanifolds}
Let $\gp$ be a simply connected, $k$-step nilpotent 
Lie group with Lie algebra $\alg.$  
Define $\qgp$ to be the simply connected, 
$(k-1)$-step nilpotent Lie group 
$\qgp=\gp/\dgp{k\hy1}.$  For $\GG$ a cocompact, 
discrete subgroup of $\gp,$
denote by $\qGG$ the image of $\GG$ 
under the canonical projection 
from $\gp$ onto $\qgp.$
The group $\qGG$ is then a cocompact, 
discrete subgroup of $\qgp.$
We denote elements of $\qgp$ by $\bar{x},$
where $\bar{x}$ is the image of $x$ under the
canonical projection from $\gp$ onto $\qgp.$
Note that $\bar{x}=\bar{y}$ if and only if $xy^{-1}$ is in 
$\dgp{k\hy1} \subset \Center.$ 
If $\Phi\in\Aut(\gp)$, we denote by $\bar\Phi\in\Aut(\qgp)$ the projection
of $\Phi$ onto $\qgp.$  
The Lie algebra $\qalg$ of $\qgp$ is just $\alg / \dalg{k\hy1}.$   
We denote elements of $\qalg$ by $\bar{X},$ 
where $\bar{X}$ is the image of $X$ under the canonical 
projection from $\alg$ onto $\qalg.$

For a left invariant metric $g$ on $\gp,$ we associate a left invariant 
metric $\qg$ on $\qgp$ by restricting $g$ to an orthogonal complement 
of $\dalg{k\hy1}$ in $\alg.$
We call the $(k\hy1)$-step nilmanifold 
$(\qNilmfld,\qg)$ the 
{\bf quotient nilmanifold} of $(\Nilmfld,g).$  
Using the definition of $\qg,$ one easily sees that the projection 
$(\Nilmfld,g)\to(\qNilmfld,\qg)$ is a Riemannian submersion
with totally geodesic fibers.

\definition{2.10 Definition} 
The simply connected nilpotent Lie group $\gp$ 
is {\bf strictly nonsingular} if for every noncentral 
$x$ in $\gp,$ if $\bar{x}=\bar{v}$ then $x$ and $v$ are conjugate
in $\gp.$  
Restated, for every noncentral $x$ in $\gp$ and every $z$ in $\Center,$
there exists $y$ in $\gp$ such that $yxy^{-1}=xz.$
The nilpotent Lie algebra $\alg$ is 
{\bf strictly nonsingular} if for every 
noncentral $X$ in $\alg$ and every $Z$ in 
$\center$ there exists $Y$ in $\alg$ such that $[X,Y]=Z.$
That is, for every noncentral $X$ in $\alg,$ 
$\center\subset\ad(X)(\alg).$
\enddefinition

\subheading{2.11 Remarks} 
\roster
\item One easily checks that these notions are equivalent, i.e., 
a nilpotent Lie group is strictly nonsingular if and only if
its Lie algebra is strictly nonsingular.
\item Equivalently, a nilpotent Lie group $\gp$ is strictly nonsingular
if and only if every fiber of the submersion
$\gp\to\qgp$ is contained in a single orbit 
of the action of $G$ on itself by conjugation.
A two-step nilpotent Lie group is strictly nonsingular if and only 
if every fiber is equal to an orbit.
\item If $\gp$ is strictly nonsingular,
then $\Center=\dgp{k\hy1},$
but not conversely.
\item The examples presented in Section 4 show that one can have
both $\gp$ and $\qgp$ 
strictly nonsingular, but this is not the case in general. 
See \cite{Gt1,2} for examples.
\endroster

\bigskip

\subheading{Section 3. New Deformations}
\nopagebreak\medskip
\nopagebreak
The method we will use to construct the new isospectral
deformations is the following. This method has previously
produced {\it discrete} families of Riemannian manifolds,
isospectral on functions but not on 1-forms \cite{G2}, \cite{Gt3,4}.

\proclaim{3.1 Theorem \cite{Gt2, Theorem 3.3}}
Let $\gp$ be a simply connected, strictly nonsingular nilpotent 
Lie group with left invariant metric $g.$  
If \ $\G1$ and $\G2$ are cocompact, discrete subgroups 
of $\gp$ such that 
$$\G1\cap\Center=\G2\cap\Center
\quad\and\quad\spec(\qnilmfld1,\qg)=
\spec(\qnilmfld2,\qg),$$
then  {$\spec(\nilmfld1,g) = \spec(\nilmfld2,g).$}
\endproclaim

\subheading{3.2 Remark} 
Any deformation using Theorem 3.1
must be of the form $(\nilmfld{s},g),$
where $\G{s}$ depends continuously on $s.$
In particular, the manifolds must have a common Riemannian covering
$(G,g)$, hence are locally isometric.
It is known \cite{Sch1} that if $\G{s}$ is any continuous 
family of cocompact, discrete subgroups of a nilpotent Lie group $\gp,$
then there exists a continuous family of automorphisms $\Phi_s\in\Aut(\gp)$
with $\Phi_0=I$ such that $\G{s}=\Phi_s(\GG).$ 
So the deformations arising from Theorem 3.1 are of the 
form $(\Nilmfld, \Phi_s^* g),$ $\Phi_s\in\Aut(\gp).$

\proclaim{3.3 Corollary} 
Let $\gp$ be a strictly nonsingular nilpotent
Lie group with left invariant metric $g$
and cocompact, discrete subgroup $\GG.$
Assume there exists a continuous family of automorphisms
$\Phi_s\in\Aut(\gp),$ $\Phi_0=I$, such that for all $s$
\roster
\item $\Phi_s$ restricted to $\Center$ is the identity, \ and 
\item $\bar\Phi_s\in\GW(\qgp, \qg; \qGG).$ 
\endroster
Then $(\Nilmfld,\Phi_s^*g)$ is an isospectral deformation.
The deformation is trivial if and only if 
$\Phi_s\in\OA(\gp,g) \circ \IA(\gp)$ for all $s.$
The deformation is Gordon--Wilson if and only if 
there exists $\Upsilon \in \OA(\gp, g)$ such that
$\Phi_s\circ\Upsilon\in\GW(\gp,g;\GG)$ for all $s.$
\endproclaim

Note again that by (2.7), if $\GW(\gp, g; \GG)$ is a group, then
we may assume $\Upsilon = I$ in the statement of Corollary 3.3.

\demo{Proof}
This follows directly from Theorem 3.1,
Theorem 2.5, and Proposition 2.8.
\qed \enddemo
\bigskip

To use Corollary 3.3 to obtain isospectral deformations that are 
not Gordon--Wilson, we exhibit a strictly nonsingular,
three-step nilpotent Lie group $\gp$ and a
continuous family of automorphisms 
$\Phi_s\in\Aut(\gp),$ $\Phi_s\not\in\GW(\gp,g)$
such that 
$\bar\Phi_s\in\GW(\qgp,\qg).$ 
For the examples presented here, 
$\AIA(\gp;\GG)= \AIA(\gp)$ and $\AIA(\qgp; \qGG)=\AIA(\qgp)$
for any choice of cocompact discrete subgroup $\GG$ of $\gp.$

\subheading{3.4 Remark}  
The key to the new deformations is close study of the behavior of the 
orthogonal automorphisms on the quotient nilmanifold.
Example I below exhibits a nontrivial isospectral deformation where 
almost inner automorphisms play no role whatsoever, 
i.e., $\GW(\qgp,\qg; \qGG)$ may be replaced by
$\OA(\qgp, \qg)$ in condition (2) of Corollary 3.3.
So the isospectral deformation $\Phi_s$ projects to a trivial
deformation on the quotient nilmanifold, but $\Phi_s$ is nontrivial.
If we disregard $\OA(\qgp,\qg),$ however, 
we reduce to the case of Gordon--Wilson deformations.
For if we let $\Phi_s$ satisfy the conditions of Corollary 3.3
and let $\bar\Phi_s\in\AIA(\qgp; \qGG)$ for all $s,$
then the elements $\bar\gamma$ and $\bar\Phi_s(\bar\gamma)$ are conjugate
for all $\bar\gamma\in\qGG.$  By strict nonsingularity, 
$\gamma$ and $\Phi_s(\gamma)$ are necessarily conjugate,
and the deformation is Gordon--Wilson.

\subheading{3.5 Remark} 
The geometric role of the orthogonal automorphisms is as follows.
Fix an orthonormal basis of the $k$-dimensional center $\center$
of $\alg$ and 
extend it to a $(k+n)$-dimensional orthonormal basis $\Cal B$ of $\alg.$ 
If $\Phi_s$ satisfies the hypotheses of Corollary 3.3,
then $\bar\Phi_s=\Upsilon_s\circ\Psi_s,$ where
$\Upsilon_s\in\OA(\qgp,\qg)$ and $\Psi_s\in\AIA(\qgp,\qGG).$
View $\Upsilon_{s*}$ as an element of $O(n)$
by considering its matrix representation with respect to $\Cal B$.  
Extend $\Upsilon_{s*}$ to $\hat\Upsilon_{s*}\in O(n+k)$ 
by the $k\times k$ identity matrix:
$$\hat\Upsilon_{s*}=
\left(\matrix \Upsilon_{s*}&0\\0&I_k\endmatrix\right).$$
The elements of the group $O(n+k),$ hence the matrices $\hat\Upsilon_{s*},$
act on our metric $g$ by taking the 
orthonormal basis $\Cal B$ of $\alg$ to $\hat\Upsilon_{s*}^{-1}\Cal B,$ 
which determines the same left invariant metric $g$ of $\gp.$
So with the obvious notation, we have
$$(\Nilmfld, \Phi_s^*g)\isom
(\Nilmfld, \hat\Psi_s^*\circ\hat\Upsilon_s^*g)\isom
(\Nilmfld, \hat\Psi_s^*g).$$
Here $\hat\Psi_{s*}=\hat\Upsilon_{s*}^{-1}\circ\Phi_s,$
a triangular matrix with 1's down the diagonal,
so $\hat\Psi_{s*}\in\SL(n+k).$ 
Note that if $\alg$ is strictly nonsingular,
then $\hat\Psi_{s*}$ and  $\hat\Upsilon_{s*}$
are Lie algebra automorphisms if and only if the 
deformation is Gordon--Wilson.

The matrices $\hat\Psi_{s*}$ deform the orthonormal basis, i.e.,
the metric $\hat\Psi_{s}^* g$ is determined by the 
orthonormal basis $\hat\Psi_{s*}^{-1}\Cal B.$
If the deformation is not Gordon--Wilson,
then the role of $\hat\Upsilon_s$ is to continuously compensate
for $\hat\Psi_s$ not being an automorphism in a way
that does not effect our metric $g,$ 
allowing us to employ Theorem 3.1.
On the other hand, assume that we know
$(\Nilmfld, \hat\Psi^*_s g)$ is an isospectral deformation
and that the  $\hat\Psi_s$ are not Lie algebra isomorphisms. 
We might errantly conclude that the manifolds 
are not locally isometric, i.e., the fact that
that the manifolds 
$(\Nilmfld, \hat\Psi^*_s g)=(\Phi_s(\GG)\backslash\gp,g)$
are covered by $(\gp,g)$  is ``hidden'' by 
the orthogonal matrices $\hat\Upsilon_{s*}.$
\bigskip

We now compute conditions for a Lie algebra automorphism of $\qalg$
to extend to a Lie group automorphism of the strictly nonsingular
Lie algebra $\alg.$ 

For the remainder, we assume the center $\center$ of $\alg$ 
is one-dimensional.
Let $\gp$ be the strictly nonsingular, nilpotent Lie group with
Lie algebra $\alg$ and left invariant metric $g.$
Let $\nn$ be the orthogonal complement of $\center$ in $\alg,$
so $\alg=\nn\oplus\center.$
Define a closed, nondegenerate two-form $\omega$ on $\nn$ by 
$$[X,Y]=([X,Y]_\nn,\omega(X,Y)\CZ)\tag{\star}$$
for all $X,Y\in\alg.$ 
Here $[X,Y]_\nn$ denotes the orthogonal
projection of $[X,Y]$ onto $\nn,$
and $\CZ$ is a unit vector in $\center.$
The condition that $\omega$ is closed is equivalent to
the Jacobi equations, and the condition that
$\omega$ is nondegenerate is equivalent to strict nonsingularity.
Note that the vector space $\nn$ with Lie bracket defined by
$[\ , \ ]_\nn$ is Lie algebra isomorphic to $\qalg.$ 

\definition{3.7 Notation}
Conversely, let $(G,g)$ be a $k$-step nilpotent Lie group.
Assume there exists a closed, nondegenerate two-form 
$\omega$ on $\alg.$  The {\bf central extension of
$\bold\alg$ by $\bold\omega$} is the strictly nonsingular, $(k+1)$-step
nilpotent Lie algebra $\halg$, where $\halg=\alg\oplus\R\CZ.$
The Lie bracket on $\halg$ is defined by $(\star),$
that is, $[(X,Z),(X',Z')]_{\halg}=([X,X']_\alg,\omega(X,X')\CZ),$
for all $(X,Z), (X',Z')$ in $\halg.$
For $X\in\alg,$ we denote by $\hX$ the element $(X,0)$ in $\halg.$
We denote by $\hgp$ the simply connected nilpotent Lie group with Lie
algebra $\halg.$
Let $g$ be a left invariant metric on $G.$
We extend $g$ to a left invariant metric $\hg$ on
$\hgp$ by requiring that $\alg$ and $\CZ$ be orthogonal
and that $\CZ$ be a unit vector. 
Note that the quotient Lie group of $(\hgp,\hg)$
is isometric to the Lie group $(G,g).$
\enddefinition

\proclaim{3.8 Proposition}
Let $\alg$ be a nilpotent Lie algebra,
and let $\omega$ be a closed, nondegenerate two-form on $\alg.$
Let $\halg$ be the central extension of $\alg$ by $\omega.$
\roster
\item A derivation
$D$ on $\alg$ extends to a derivation $\hat{D}$
on $\halg$ if and only if there exists $\eta\in\alg^*$ such that
$D^*\omega=d\eta,$ i.e., $D^*\omega$ is exact.
\item A skew-symmetric derivation
$S$ on $\alg$ extends to a skew-symmetric derivation $\hat{S}$
on $\halg$ if and only if $S^*\omega=0.$
\endroster
\endproclaim

\definition{3.9 Notation} If $D^*\omega=d\eta$ with $\eta\in\alg^*,$ 
we denote by $\hat{D}$ the derivation of
$\halg=\alg\oplus\R\CZ$
defined by  
$\hat{D}(\hX)=(D(X),-\eta(X)\CZ)$ for all $X\in\alg,$
and $\hat{D}(\CZ)=0.$ 
\enddefinition

\demo{Proof of 3.8}
A straightforward calculation using the definition of
a derivation shows that $D$ extends to a derivation $\hat{D}$
on $\halg$ if and only if $\hat{D}$ is as in (3.9), and (1) follows.
The proof of (2) uses (1) and the definition of $\hg.$
\qed \enddemo

\proclaim{3.10 Corollary} 
Let $\gp$ be a strictly nonsingular nilpotent
Lie group with Lie algebra $\alg$ and left invariant metric $g.$
Assume there exists a closed, nondegenerate two-form
$\omega$ on $\alg$ and a derivation $D\in\Der(\alg)$ such that 
\roster
\item $D=S+A$ where $S\in \OD(\alg,g)$ and $A\in\AID(\alg),$ and
\item $D^*\omega$ is exact.
\endroster
Denote by $\hat{D}\in\Der(\halg)$ the extension of $D$ to $\halg.$
Define the automorphisms
$\Phi_s\in\Aut(\hgp)$ on the Lie algebra level
by $\Phi_{s*}=\exp(s\hat{D}).$ 
Let $\hGG$ be any cocompact, discrete subgroup of $\hgp.$
Then $(\hNilmfld,\Phi_s^*\hg)$ is an isospectral 
deformation. 
The deformation is trivial if and only if 
$S^*\omega=0$ and $A\in\ID(\alg).$
The deformation is Gordon--Wilson 
if and only if $S^*\omega=0.$
\endproclaim

\demo{Proof of 3.10}
Follows from Corollary 3.3, Proposition 3.8, and Remark 2.7.
Note that by strict nonsingularity,
if $\Psi\in\Aut(\gp)$ and $\bar\Psi\in\AIA(\qgp),$
then $\Psi\in\AIA(\gp).$
\qed \enddemo

\subheading{3.11 Remark}
With appropriate additional assumptions,
one may extend the above Corollary to 
include $\AID(\alg; \GG)$. The statement becomes somewhat cumbersome,
however, and so is not included here.
\bigskip

\subheading{Section 4: Examples}
\medskip

Consider the strictly nonsingular, two-step nilpotent Lie algebra 
$\alg$ with orthonormal basis $\{X_1,$ $ X_2,$ $ X_3,$ $ X_4,$ $ Z_1,$ $Z_2\}$ 
and Lie bracket given by
$$\align
[X_1,X_2] = [X_3,X_4] &= Z_1,\\
[X_1,X_3] = [X_4,X_2] &= Z_2,
\endalign$$
and all other basis brackets zero. 
Let $\{\a_1, \a_2, \a_3, \a_4, \z_1, \z_2\}$ be the dual basis
of $\alg^*,$
and let $\gp$ be the simply
connected Lie group with Lie algebra $\alg.$

\subheading{Example I}

Let $$\omega = \a_1\wedge\a_2+\a_2\wedge\a_3+
\a_3\wedge\a_4+\a_1\wedge\z_1-\a_4\wedge\z_2.$$
The closed, nondegenerate two-form $\omega$ 
induces a strictly nonsingular, three-step 
nilpotent Lie algebra $\halg=\alg\oplus\R\CZ$ 
with Lie bracket and left invariant metric as defined
in (3.7). So $\halg$
is the 7-dimensional Lie algebra with orthonormal basis 
$\{\hX_1, \hX_2, \hX_3, \hX_4, \hZ_1, \hZ_2, \CZ\}$
and Lie bracket
$$\align
[\hX_1,\hX_2]&=[\hX_3,\hX_4]=\hZ_1+\CZ,\\
[\hX_1,\hX_3]&=[\hX_4,\hX_2]=\hZ_2,\\
[\hX_2,\hX_3]&=[\hX_1,\hZ_1]=[\hZ_2,\hX_4]=\CZ,
\endalign$$
and all other basis brackets zero.

One easily checks that 
the skew-symmetric derivation $S:\alg\to\alg$ defined by 
$$X_1\to -X_4;\quad X_2\to 2 X_3;
\quad X_3\to-2 X_2;\quad X_4\to X_1;\quad
Z_1\to Z_2;\quad Z_2\to -Z_1$$
satisfies $S^*\omega=d\z_2,$ so $S$ extends to the derivation
$\hat{S}:\halg\to\halg$ given by 
$\hat{S}(\hX)=(S(X),-\z_2(X)\CZ)$, $\hat{S}(\CZ)=0.$  That is, 
$\hat{S}$ sends
$$\hX_1\to -\hX_4;\quad\hX_2\to 2\hX_3;
\quad\hX_3\to-2\hX_2;\quad\hX_4\to \hX_1; \quad
\hZ_1\to\hZ_2;\quad \hZ_2\to -\hZ_1-\CZ; \quad \CZ\to0.$$

Let $\hgp$ denote the simply connected, strictly nonsingular nilpotent Lie
group with Lie algebra $\halg.$
Define a continuous family $\Phi_s$ of Lie group automorphisms
on $\hgp$ given at the Lie algebra level by
$\Phi_{s*}=\exp(s\hat{S}).$
Using any symbolic algebra package, one computes
that the automorphism $\Phi_{s*}$ sends
$$\align
\hX_1&\to \cos(s)\hX_1-\sin(s)\hX_4,\\
\hX_2&\to \cos(2s)\hX_2+\sin(2s)\hX_3,\\
\hX_3&\to -\sin(2s)\hX_2+\cos(2s)\hX_3,\\
\hX_4&\to \sin(s)\hX_1+\cos(s)\hX_4,\\
\hZ_1&\to \cos(s)\hZ_1+\sin(s)\hZ_2-(1-\cos(s))\CZ,\\
\hZ_2&\to -\sin(s)\hZ_1+\cos(s)\hZ_2-\sin(s)\CZ,\\
\CZ&\to\CZ.\endalign$$

Let $\hGG$ be any cocompact, discrete subgroup of $\hgp$
and let $\hG{s}=\Phi_s(\hGG).$
Then by Corollary 3.10,
$(\hnilmfld{s},\hg)=(\hNilmfld,\Phi_s^*\hg)$
is an isospectral deformation that is not Gordon--Wilson.
Theorem 6.1 below shows that  
$1\hy\spec(\Phi_s(\hGG)\backslash \hgp,\hg)$
is nonconstant, providing an alternate proof that
the deformation is not through almost inner automorphisms.

By factoring out the
orthogonal matrices as in (3.5), 
the left invariant metric $\Phi_s^*\hg$ is determined by the orthonormal basis
$\{\hX_1, \hX_2, \hX_3, \hX_4, \hZ_1+(1-\cos(s))\CZ, 
\hZ_2+\sin(s)\CZ, \CZ \}$ of $\halg.$
\bigskip

\subheading{Example II}

Let $$\omega=\a_2\wedge\a_3+\a_2\wedge\a_4+\a_1\wedge\z_1-\a_4\wedge\z_2.$$

The closed, nondegenerate two-form $\omega$ 
induces a strictly nonsingular, three-step 
nilpotent Lie algebra $\halg=\alg\oplus\R\CZ,$ 
with Lie bracket and left invariant metric as defined
in (3.7).  So $\alg$ is the 7-dimensional Lie algebra with
orthonormal basis
$\{\hX_1, \hX_2, \hX_3, \hX_4, \hZ_1, \hZ_2, \CZ\},$
and Lie bracket
$$\align
[\hX_1,\hX_2]&=[\hX_3,\hX_4]=\hZ_1,\\
[\hX_1,\hX_3]&=\hZ_2, \ [\hX_4,\hX_2]=\hZ_2-\CZ,\\
[\hX_2,\hX_3]&=[\hX_1,\hZ_1]=[\hZ_2,\hX_4]=\CZ,
\endalign$$
and all other basis brackets zero.

One easily checks that the derivations
$S:\alg\to\alg$ defined by 
$$X_1\to -X_4;\quad X_2\to 2X_3;
\quad X_3\to -2X_2;\quad X_4\to X_1; \quad Z_1\to Z_2;\quad Z_2\to -Z_1;$$
and $A:\alg\to\alg$ defined by 
$$X_1, X_4, Z_1, Z_2 \to 0; \quad X_2\to -Z_1; \quad X_3\to 2 Z_2 $$
satisfy $S^*\omega=\a_1\wedge\a_2-2\a_3\wedge\a_4,$
and $A^*\omega=-S^*\omega.$
Note that $S^*\omega$ is 
not in $\span_\R\{d\z_1,d\z_2,d\xi\}$, hence not exact.
So neither $S$ nor $A$ will extend to a derivation on 
$\halg.$ However, the sum $D=A+S$ is exact as $D^*\omega=0,$
so $D$ extends to the derivation $\hat{D}:\halg\to\halg$
given by $\hat{D}(\hX)=(S(X)+A(X),0)$, $\hat{D}(\CZ)=0.$
That is, $\hat{D}$ sends
$$\hX_1\to -\hX_4;\quad\hX_2\to 2\hX_3-\hZ_1;
\quad\hX_3\to-2\hX_2+ 2 \hZ_2;\quad\hX_4\to \hX_1;$$
$$\hZ_1\to\hZ_2;\quad \hZ_2\to -\hZ_1; \quad \CZ\to0.$$

Let $\hgp$ denote the simply connected, strictly nonsingular nilpotent Lie
group with Lie algebra $\halg.$
Define a continuous family $\Phi_s$ of Lie group automorphisms
of $\hgp$ defined on the Lie algebra level by
$\Phi_{s*}=\exp(s\hat{D}).$
The automorphism $\Phi_{s*}$ sends
$$\align
\hX_1&\to \cos(s)\hX_1-\sin(s)\hX_4,\\
\hX_2&\to \cos(2s)\hX_2+\sin(2s)\hX_3-\sin(s)\hZ_1+(\cos(s)-\cos(2s))\hZ_2,\\
\hX_3&\to -\sin(2s)\hX_2+\cos(2s)\hX_3+\sin(2s)\hZ_2,\\
\hX_4&\to \sin(s)\hX_1+\cos(s)\hX_4,\\
\hZ_1&\to \cos(s)\hZ_1+\sin(s)\hZ_2,\\
\hZ_2&\to -\sin(s)\hZ_1+\cos(s)\hZ_2,\\
\CZ&\to\CZ.\endalign$$

Let $\hGG$ be any cocompact, discrete subgroup of $\hgp.$
Then 
$(\Phi_s(\hGG)\backslash \hgp,\hg)=(\hNilmfld, \Phi_s^*\hg)$
is an isospectral deformation not arising from almost inner automorphisms.
We show in Theorem 6.2 below that  
$1\hy\spec(\Phi_s(\hGG)\backslash \hgp,\hg)$
is nonconstant, providing an alternate proof that
the deformation is not Gordon--Wilson.

By factoring out the
orthogonal matrices as in (3.5), 
we see that $\Phi_s^*\hg$ is determined by the orthonormal basis
$\{\hX_1,$ 
$\hX_2+\sin(s)\cos(2s)\hZ_1-(1-\cos(s)\cos(2s))\hZ_2,$ 
$\hX_3-\sin(s)\sin(2s)\hZ_1-\cos(s)\sin(2s)\hZ_2,$ 
$\hX_4,$ $\hZ_1,$ $\hZ_2,$ $\CZ \}$ of $\halg.$

\bigskip

\subheading{Section 5.  The Marked Length Spectrum}
\medskip

In this section we prove the following.

\proclaim{5.1 Theorem}
The nilmanifolds $(\hnilmfld{s},\hg)$ presented in Example I above
have the same marked length spectrum for all $s.$
In particular, $\Phi_s$ marks the length spectrum from
$(\hnilmfld0,\hg)$ to $(\hnilmfld{s},\hg).$ 
\endproclaim

\proclaim{5.2 Theorem}
The nilmanifolds $(\hnilmfld{s},\hg)$ presented in Example II above
have the same marked length spectrum for all $s.$
In particular, $\Phi_s$ marks the length spectrum from
$(\hnilmfld0,\hg)$ to $(\hnilmfld{s},\hg).$ 
\endproclaim

\subheading{5.3 Remark}
These are the first continuous
families of Riemannian manifolds with the same marked length
spectrum but not the same spectrum on 1-forms. 
Discrete families of manifolds with these properties were presented
in \cite{Gt4}. The only other known examples of manifolds that have the 
same marked length spectrum but not the same 1-form spectrum
are the Zoll spheres and standard spheres, which are 
not even isospectral on functions.
\smallskip

If $\gp$ is a simply connected Lie group with 
cocompact, discrete subgroup $\GG,$ then $\GG$ is also the fundamental
group of $\Nilmfld$,
so free homotopy classes of loops of $\Nilmfld$ correspond
to conjugacy classes in $\GG.$ 
We denote by $[\gamma]_\GG$ the conjugacy class in $\GG$ represented
by $\gamma\in\GG,$ that is, 
$[\gamma]_\GG=\{\gamma'\gamma{\gamma'}^{-1}:\gamma'\in\GG\}.$
Recall that any isomorphism between fundamental groups 
induces a correspondence between free homotopy classes of loops.

\definition{5.4 Definition} A pair of Riemannian manifolds 
$(\G1\backslash \tilde{M}_1, g_1)$ and 
$(\G2\backslash \tilde{M}_2, g_2)$ have the same {\bf marked length spectrum}
if and only if there exists an isomorphism $\Phi:\G1 \rightarrow\G2$
between their fundamental groups such that
corresponding free homotopy classes contain closed geodesics
of the same length. That is, 
the free homotopy class represented by 
$[\gamma]_{\G1}$ contains a closed geodesic of
(positive) period $\lambda$ if and only if
the free homotopy class represented by 
$[\Phi(\gamma)]_{\G2}$ also contains a closed geodesic of
period $\lambda$. We say that $\Phi$ {\bf marks the length spectrum}
from $(\G1\backslash\tilde{M}_1, g_1)$ to 
$(\G2\backslash\tilde{M}_2, g_2).$
\enddefinition

We have the following sufficient condition 
for a pair of nilmanifolds to have the same marked length spectrum.

\proclaim{5.5 Theorem \cite{G1, Theorem 1.8}}
Let $\gp$ be a simply connected nilpotent Lie group with cocompact,
discrete subgroup $\GG.$
Let $\Phi$ be a $\GG$-almost inner automorphism of $\gp.$
Then for any left invariant metric $g$ of $G,$
$\Phi$ marks the length spectrum from $(\Nilmfld,g)$
to $(\Phi(\GG)\backslash\gp,g).$
\endproclaim

Trivially, if $\Phi$ is an orthogonal automorphism, i.e., $\Phi\in\OA(\gp,g),$
then $\Phi$ marks the length spectrum from
$(\Nilmfld,g)$ to $(\Phi(\GG)\backslash\gp,g).$
Conversely, on two-step nilmanifolds Eberlein showed the following.

\proclaim{5.6 Theorem \cite {E, Theorem 5.20}}
Let $\G1,\G2$ be cocompact, discrete subgroups of simply connected,
two-step nilpotent Lie groups $\gp_1$, $\gp_2$ with left invariant metrics 
$g_1$ , $g_2$ respectively.
Assume that  $(\G1\backslash\gp_1,g_1)$ and 
$(\G2\backslash \gp_2,g_2)$ 
have the same marked length spectrum, and
let $\Phi:\G1\rightarrow\G2$ be an isomorphism inducing
this marking.
Then $\Phi = {(\Phi_1\circ\Phi_2)\vert}_{\G1},$ 
where $\Phi_2$ is a $\G1$-almost inner automorphism of $\gp_1$,
and $\Phi_1$ is an isomorphism of $(\gp_1,g_1)$ onto $(\gp_2,g_2)$
that is also an isometry. Moreover, this factorization is unique.
In particular, if $(\G1\backslash\gp_1,g_1)$ and
$(\G2\backslash\gp_2,g_2)$ have the same marked length spectrum, 
they necessarily have the same spectrum of the Laplacian on functions
and on $p$-forms for all $p.$
\endproclaim

In contrast to (5.6), we have the following necessary and 
sufficient condition for a large family of pairs of three-step 
nilmanifolds to have the same marked length spectrum.

\proclaim{5.7 Theorem \cite{Gt4, Theorem 3.3.1}} 
Let $\gp$ be a simply connected, 
strictly nonsingular, three-step nilpotent Lie group 
with a one-dimensional center.  Let $\G1$ and $\G2$ 
be cocompact, discrete subgroups of $\gp$ such that 
$\G1\cap\Center=\G2\cap\Center.$ Let $g$ be any 
left invariant metric on $\gp.$ Then 
$(\nilmfld1,g)$ and $(\nilmfld2,g)$ 
have the same marked length spectrum 
if and only if 
there exists an isomorphism 
$\Phi:\G1\to\G2$ such that 
$\bar\Phi:\qG1\to\qG2$ 
marks the length spectrum between 
$(\qnilmfld1,\qg)$ and $(\qnilmfld2,\qg).$ 
\endproclaim

\demo{Proof of 5.1 and 5.2}
The automorphisms $\Phi_s$ presented in Example I and Example
II have the property that $\bar\Phi_s$ is the composition of
an almost inner automorphism and an orthogonal automorphism,
which by (5.5) necessarily mark the length spectrum from
$(\nilmfld0, g)$ to $(\nilmfld{s}, g).$ 
The result now follows from (5.7).
\qed \enddemo

We also showed in \cite{Gt4} that for a large class of three-step nilmanifolds,
which include Example I and Example II,
if a family of nilmanifolds in this class have the same marked length
spectrum, they necessarily share the same Laplace spectrum on functions,
providing an alternate proof that Example I and Example II
are isospectral deformations.

\proclaim{5.8 Theorem \cite{Gt4, Theorem 3.3.2}}
Let $\gp$ be a simply connected, strictly nonsingular, 
three-step nilpotent Lie group.  Let $\G1$ and $\G2$ be cocompact, 
discrete
subgroups of $\gp$  such that $\G1\cap \Center = \G2\cap \Center.$ 
If $(\nilmfld1, g)$ and $(\nilmfld2, g)$ 
have the same marked length spectrum, 
then $(\nilmfld1, g)$ and $(\nilmfld2, g)$ 
are isospectral on functions.
\endproclaim

\bigskip

\subheading{Section 6:  Comparing the 1-form spectrum of nilmanifolds}
\medskip

In this section, we prove the following. 

\proclaim{6.1 Theorem}
The nilmanifolds $(\hnilmfld{s}, \hg)$ of Example I  
are not isospectral on 1-forms.  
\endproclaim
\medskip

\proclaim{6.2 Theorem}
The nilmanifolds $(\hnilmfld{s}, \hg)$ of Example II
are not isospectral on 1-forms.  
\endproclaim
\medskip

\demo{Proof of Theorem 6.1}

Let $\hgp,$ $\halg,$ $\hg,$ and $\Phi_s$ be as in Example I, Section 4.
Let $\hGG$ be an arbitrary cocompact, discrete subgroup of $\hgp,$
and let $\hG{s}=\Phi_s(\hGG).$
We calculate a monic polynomial $P_s(x)$ 
of degree 7 whose roots must appear
in $1\hy\spec(\hnilmfld{s}, \hg)$.  We show that the
coefficients of $P_s(x)$, which must change continuously
with $s,$ are nonconstant. The result follows.

We refer the reader to \cite{Gt3, Appendix} for 
definitions, notation, and an exposition of the
method used below to demonstrate the inequality of the 1-form spectrum.
We include only those details needed for verification. 

The main tool we use is representation theory.
Irreducible representations of nilpotent Lie groups are
well understood due to Kirillov theory.
In particular, the irreducible representations of a nilpotent
Lie group $\gp$ correspond to the coadjoint orbits of the
Lie algebra $\alg$ of $\gp.$
For $\tau\in\alg^*,$ we denote by $\pi_\tau$ the corresponding
irreducible representation of $\gp$.
A good reference for representations of nilpotent Lie groups 
is \cite{CG}.

We use only the characters of $\hgp,$
i.e., the (complex) one-dimensional irreducible representations
to demonstrate the inequality of the 1-form spectrum.
The dimension of the representation space $\H_\tau$ of $\pi_\tau$ 
corresponds to the dimension of the coadjoint orbit of $\tau.$ 
In particular, the one-dimensional representations of $\gp$
correspond to those elements of $\alg^*$ that are zero
on the derived algebra $[\alg, \alg]$ of $\alg.$

Consider $\hrep{s},$ the quasi-regular representation 
of $\hgp$ on $\hfnsp{s}.$
For $f \in \hfnsp{s},$ and $x\in\hgp,$  $\hrep{s}(x)f=f\circ R_x,$
where $R$ denotes the right action of $\hgp$ on $\hnilmfld{s}.$
It is known that $\hfnsp{s}$ is completely reducible as
$$\hfnsp{s}\cong{\underset{\tau\in\Tau_{s}}\to\bigoplus}
{m_s(\tau)}\H_\tau.$$
Here $m_s(\tau)$ is the multiplicity of $\pi_\tau$ in $\hrep{s},$
and $\Tau_s$ is a subset of $\halg^*$ containing at most one element 
of each coadjoint orbit of $\halg^*.$ 

Let $\{\a_1,$  $\a_2,$ $\a_3,$ $\a_4,$ $\z_1,$ $\z_2,$ $\xi \}$ 
be the dual basis to the orthonormal basis
$\{\hX_1,$ $\hX_2,$ $\hX_3,$ $\hX_4,$ $\hZ_1,$ $\hZ_2,$ $\CZ \}$ of $\halg.$
Denote by $\Char$ the subspace $\Char=\span_\R\{\a_1, \a_2, \a_3, \a_4\}$
of $\halg^*.$
Note that $\Char=\{\tau\in\halg^*:\ \tau([\halg,\halg]) \equiv 0\}.$
Thus the characters of $\hgp$ appearing in $\hfnsp{s}$ 
correspond to $\Tau_s \cap \Char.$

Let $\tau\in\Char.$
We now calculate $\Delta$ acting on 
$\H_\tau \otimes \Lambda^1(\halg^*),$
an invariant subspace of $\Delta.$
Note that if $\Cplx \Lambda^1(\halg^*)$ denotes $\Lambda^1(\halg^*)$ 
with complex coefficients, then 
$\H_\tau \otimes \Lambda^1(\halg^*) = 
\F \otimes \Cplx \Lambda^1(\halg^*),$ where
$\F(x)= e^{(2\pi i \tau(\log(x)))}.$

As $\hZ_1\F=\hZ_2\F=\CZ\F=0,$
we know $\Delta\F= 4\pi^2||\tau||^2\F.$
One also computes the following.
$$\Delta\a_1 = \Delta\a_2 = \Delta\a_3 = \Delta\a_4 =0, \quad
\Delta \z_1 = 2\z_1 + 2 \xi, \quad \Delta\z_2 = 2 \z_2, \quad
\Delta \xi = 2 \z_1 + 5\xi$$
{\eightpoint
$$
\vbox{\offinterlineskip
\halign{\strut\vrule#&\quad#\hfil\quad&&\vrule#&\quad\hfil#\hfil\quad\cr
\noalign{\hrule}
&&&&&&&&&&&&&&&&\cr
&$\nabla_U\mu$&&$\a_1$
&&$\a_2$&&$\a_3$&&$\a_4$&&$\z_1$&&$\z_2$&&$\xi$&\cr
&&&&&&&&&&&&&&&&\cr
\noalign{\hrule}
\noalign{\hrule}
&&&&&&&&&&&&&&&&\cr
&$X_1$
&&$0$
&&$\frac1{2}(\z_1+\xi)$
&&$\frac1{2}\z_2$
&&$0$
&&$\frac1{2}(-\a_2+\xi)$
&&$-\frac1{2}\a_3$
&&$-\frac1{2}(\a_2+\z_1)$
&\cr
&&&&&&&&&&&&&&&&\cr
\noalign{\hrule}
&&&&&&&&&&&&&&&&\cr
&$X_2$
&&$-\frac1{2}(\z_1+\xi)$
&&$0$
&&$\frac1{2}\xi$
&&$-\frac1{2}\z_2$
&&$\frac1{2}\a_1$
&&$\frac1{2}\a_4$
&&$\frac1{2}(\a_1-\a_3)$
&\cr
&&&&&&&&&&&&&&&&\cr
\noalign{\hrule}
&&&&&&&&&&&&&&&&\cr
&$X_3$
&&$-\frac1{2}\z_2$
&&$-\frac1{2}\xi$
&&$0$
&&$\frac1{2}(\z_1+\xi)$
&&$-\frac1{2}\a_4$
&&$\frac1{2}\a_1$
&&$\frac1{2}(\a_2-\a_4)$
&\cr
&&&&&&&&&&&&&&&&\cr
\noalign{\hrule}
&&&&&&&&&&&&&&&&\cr
&$X_4$
&&$0$
&&$\frac1{2}\z_2$
&&$-\frac1{2}(\z_1+\xi)$
&&$0$
&&$\frac1{2}\a_3$
&&$-\frac1{2}(\a_2+\xi)$
&&$\frac1{2}(\a_3+\z_2)$
&\cr
&&&&&&&&&&&&&&&&\cr
\noalign{\hrule}
}}
$$
}
\medskip 

Write $\tau=A_1\a_1+A_2\a_2+A_3\a_3+A_4\a_4$
where $A_i\in\R.$
Using $\Delta(f\tau)=(\Delta f)\tau + f(\Delta\tau) -2\nabla_{\grad f} \tau$
and the above information, 
a straightforward calculation shows that if we set
$E_\tau = $
$$\left(
\matrix 
4\pi^2||\tau||^2 & 0 & 0 & 0 &  2\pi i A_2 & 2\pi i A_3 & 2\pi i A_2\\
0& 4 \pi^2||\tau||^2 & 0 & 0 & -2\pi i A_1 & -2\pi i A_4 & 2\pi i(A_3-A_1)\\
0&0& 4 \pi^2||\tau||^2 & 0 & 2\pi i A_4 & -2\pi i A_1 & 2\pi i(A_4-A_2)\\
0&0&0& 4 \pi^2||\tau||^2 & -2 \pi i A_3& 2 \pi i A_2 & -2\pi i A_3\\
-2\pi i A_2&2\pi i A_1&-2\pi i A_4&2\pi i A_3&  4 \pi ^2 ||\tau||^2 + 2 &0 & -2 \pi i A_1 +2\\
-2\pi i A_3&2\pi i A_4&2\pi i A_1&-2\pi i A_2&0&  4 \pi ^2 ||\tau||^2 + 2 & 2 \pi i A_4 \\
-2\pi i A_2&2\pi i (A_1-A_3)&2\pi i (A_2-A_4)&2\pi i A_3&2\pi i A_1+2&-2\pi i A_4& 4 \pi^2||\tau||^2 + 5
\endmatrix
\right)$$
then there exists
$\mu\in\Cplx \Lambda^1(\halg^*),$ such that
$\Delta(\F \otimes \mu) = \lambda(\F \otimes \mu)$ if and only if 
$\lambda$ is an eigenvalue of the matrix  $E_\tau.$
That is, 
for every $\tau\in\Tau_s\cap\Char,$ 
the roots of the characteristic polynomial of $E_\tau$
are contained in $1\hy\spec(\nilmfld{s},g).$
One easily computes that 
$\tau\in\Tau_s \cap\Char$ if and only if 
$\tau(\log \hG{s}) \subset \Z.$ 

Now $\log\hG{s}$ generates a lattice of full rank in $\halg,$ so
we can find $V^s_1,$ $V^s_2,$ $V^s_3,$ $V^s_4,$ elements of 
$\span_\R\{\hX_1,$ $\hX_2,$ $\hX_3,$ $\hX_4\},$ with the property that
$\tau(\log\hG{s})\subset\Z$
if and only if $\tau(V^s_i)\subset \Z,$ $i=1,2,3,4.$
Let $L_s=\span_\Z\{V^s_1,V^s_2,V^s_3,V^s_4\}.$
Note that $L_s$ is a lattice of full rank in 
$\span_\R\{\hX_1,\hX_2,\hX_3,\hX_4\}.$
Let $L_s^*$ be the dual lattice of $L_s$
in $\Char,$ so $L_s^*$ is a lattice of full rank in $\Char,$ and
$$\tau\in\Tau_s\cap\Char \ \text{ if and only if } \ \tau \in L_s^*.$$

A straightforward computation shows that if
$\tau=a_1\a_1+a_2\a_2+a_3\a_3+a_4\a_4$
is any element of $L_0^*\subset\Tau_0\cap\Char$ then 
$\tau_s=\tau\circ \Phi^{-1}_{s*}$ 
is an element of $L_s^*\subset\Tau_s\cap\Char.$ 
That is, 
$\tau_s=A_1(s)\a_1+A_2(s)\a_2+A_3(s)\a_3+A_4(s)\a_4$
is an element of $L_s^*,$ where
$$
\align
A_1(s) &=a_1\cos(s)+a_4\sin(s)\\
A_2(s) &=a_2\cos(2s)-a_3\sin(2s)\\
A_3(s) &=a_2\sin(2s)+a_3\cos(2s)\\
A_4(s) &=-a_1\sin(s)+a_4\cos(s).
\endalign$$

Note that $||\tau_s||^2=a_1^2+a_2^2+a_3^2+a_4^2$ is independent of $s.$
Using any symbolic computation package, one easily checks
that the characteristic polynomial of $E_{\tau_s},$ denoted $P_s(x),$
is a monic polynomial of degree 7, and the coefficients
of the monic polynomial $P_s(x)$ are
independent of $s$ if and only if $a_1 a_2=a_3 a_4.$
As $L_0^*$ is full lattice in $\Char,$ there exists
at least one $\tau\in L_0^*$ satisfying $a_1 a_2 \neq a_3 a_4,$
and we are done.
\qed \enddemo

\subheading{6.3 Remark}
Alternatively, one may perform the same computation by viewing
the manifolds of Example I as $(\hNilmfld, \Phi_s^*\hg).$
This is simplified by factoring out an orthogonal matrix for each $s,$
as in (3.5), obtaining the orthonormal basis 
$\{\hX_1, \hX_2, \hX_3, \hX_4, \hZ_1+(1-\cos(s))\CZ, 
\hZ_2+\sin(s)\CZ, \CZ \}$ of $\halg$.
Here the functional $\tau$ remains fixed, 
but the values for $\Delta \a_i$, $\Delta \z_i$ and $\nabla_U\mu$
listed above now vary with $s.$ 
One obtains a different expression for $E_{\tau}$ depending on $s,$
but the characteristic polynomial $P_s(x)$ is equal to that obtained 
above. 

\demo{Proof of Theorem 6.2}

Let $\hgp,$ $\halg,$ $\hg,$ and $\Phi_s$ be as in Example II, Section 4.
Let $\hGG$ be an arbitrary cocompact, discrete subgroup of $\hgp,$
and let $\hG{s}=\Phi_s(\hGG).$
Let $\{\a_1, \a_2, \a_3, \a_4, \z_1, \z_2, \xi \}$ 
be the dual basis to the orthonormal basis
$\{\hX_1, \hX_2, \hX_3, \hX_4, \hZ_1, \hZ_2, \CZ \}$ of $\halg.$

The proof here is identical to that of Example I, with the following
modifications.
$$\Delta\a_1 = \Delta\a_2 = \Delta\a_3 = \Delta\a_4 =0,
\quad  \Delta \z_1 = 2\z_1, \quad  \Delta\z_2 = 2 \z_2 - \xi, \quad 
\Delta \xi = -\z_1 + 4\xi.$$
{\eightpoint
$$
\vbox{\offinterlineskip
\halign{\strut\vrule#&\quad#\hfil\quad&&\vrule#&\quad\hfil#\hfil\quad\cr
\noalign{\hrule}
&&&&&&&&&&&&&&&&\cr
&$\nabla_U\mu$&&$\a_1$
&&$\a_2$&&$\a_3$&&$\a_4$&&$\z_1$&&$\z_2$&&$\xi$&\cr
&&&&&&&&&&&&&&&&\cr
\noalign{\hrule}
\noalign{\hrule}
&&&&&&&&&&&&&&&&\cr
&$X_1$
&&$0$
&&$\frac1{2}\z_1$
&&$\frac1{2}\z_2$
&&$0$
&&$\frac1{2}(-\a_2+\xi)$
&&$-\frac1{2}\a_3$
&&$-\frac1{2}\z_1$
&\cr
&&&&&&&&&&&&&&&&\cr
\noalign{\hrule}
&&&&&&&&&&&&&&&&\cr
&$X_2$
&&$-\frac1{2}\z_1$
&&$0$
&&$\frac1{2}\xi$
&&$\frac1{2}(-\z_2+\xi)$
&&$\frac1{2}\a_1$
&&$\frac1{2}\a_4$
&&$-\frac1{2}(\a_3+\a_4)$
&\cr
&&&&&&&&&&&&&&&&\cr
\noalign{\hrule}
&&&&&&&&&&&&&&&&\cr
&$X_3$
&&$-\frac1{2}\z_2$
&&$-\frac1{2}\xi$
&&$0$
&&$\frac1{2}\z_1$
&&$-\frac1{2}\a_4$
&&$\frac1{2}\a_1$
&&$\frac1{2}\a_2$
&\cr
&&&&&&&&&&&&&&&&\cr
\noalign{\hrule}
&&&&&&&&&&&&&&&&\cr
&$X_4$
&&$0$
&&$\frac1{2}(\z_2-\xi)$
&&$-\frac1{2}\z_1$
&&$0$
&&$\frac1{2}\a_3$
&&$-\frac1{2}(\a_2+\xi)$
&&$\frac1{2}(\a_2+\z_2)$
&\cr
&&&&&&&&&&&&&&&&\cr
\noalign{\hrule}
}}
$$
}
Set $E_\tau=$
$$\left(
\matrix 
4\pi^2||\tau||^2 & 0 & 0 & 0 &  2\pi i A_2 & 2\pi i A_3 & 0\\
0& 4 \pi^2||\tau||^2 & 0 & 0 & -2\pi i A_1 & -2\pi i A_4 & 2\pi i(A_3+A_4)\\
0&0& 4 \pi^2||\tau||^2 & 0 & 2\pi i A_4 & -2\pi i A_1 & -2\pi i A_2\\
0&0&0& 4 \pi^2||\tau||^2 & -2 \pi i A_3& 2 \pi i A_2 & -2\pi i A_2\\
-2\pi i A_2&2\pi i A_1&-2\pi i A_4&2\pi i A_3&  4 \pi ^2 ||\tau||^2 + 2 &0 & -2 \pi i A_1\\
-2\pi i A_3&2\pi i A_4&2\pi i A_1&-2\pi i A_2&0&  4 \pi ^2 ||\tau||^2 + 2 & 2 \pi i A_4 - 1 \\
0&-2\pi i (A_3+A_4)&2\pi i A_2&2\pi i A_2&2\pi i A_1&-2\pi i A_4-1& 4 \pi^2||\tau||^2 + 4
\endmatrix
\right)$$
We may use the same values for $A_i(s),$ $i=1,2,3,4.$
\qed \enddemo

\enddocument